\documentclass[10pt]{article}

\usepackage[a4paper,margin=1in]{geometry}
\usepackage{setspace}
\usepackage{amsmath,amssymb,amsthm,mathtools}

\usepackage{graphicx}
\usepackage{subcaption}

\usepackage{algorithm}
\usepackage{algpseudocode}

\usepackage{booktabs}

\usepackage{xcolor}
\usepackage{enumitem}
\usepackage{csquotes}

\usepackage[hidelinks]{hyperref}
\usepackage{url}
\usepackage[numbers,sort&compress]{natbib}
\title{A Computationally Tractable Path-Planning Method for Airborne Wind Energy Systems}

\author{
Manuel C.R.M. Fernandes\thanks{SYSTEC, Faculdade de Engenharia, Universidade do Porto, Porto, Portugal. \texttt{mcrmf@fe.up.pt}}
\and
Fernando A.C.C. Fontes\thanks{SYSTEC, Faculdade de Engenharia, Universidade do Porto, Porto, Portugal. \texttt{faf@fe.up.pt}}
}

\date{\today}

\begin{document}

\maketitle

\begin{abstract}
Airborne Wind Energy Systems (AWES) have emerged as a promising renewable energy technology that exploits stronger, more consistent high-altitude winds via tethered airborne devices. Among the various concepts, crosswind systems, where efficient flight control is essential to maximise energy output, offer significant potential. This paper addresses the problem of reference selection for crosswind flight control, focusing on the design of power-maximising geometric flight paths for the reel-out phase of Groundgen systems.

To overcome the computational challenges associated with optimal control approaches, a computationally tractable framework is proposed in which a path-planning problem is formulated as a nonlinear program. The method optimises the parameters of a Lissajous curve to maximise the average power production over the reel-out phase, while incorporating curvature constraints. The proposed approach provides an efficient alternative to existing optimal control and learning-based methods.

\end{abstract}

\section{Introduction}

Airborne Wind Energy Systems (AWES) are autonomous tethered flying devices that aim to harvest high-altitude winds for power generation. These systems, due to their low material and structural requirements, have been proposed as a potentially game-changing technology in renewable energy systems, notably for their greater transportability and their capability of achieving winds at higher altitudes than conventional wind turbines, where the wind speeds are stronger and more consistent. For a general view of the technology, its challenges and current stage of development, we point to the works  \cite{Cherubini_2015}, \cite{Vermillion_2021}, and \cite{EuropeanCommission_2018}. 

Although there are several different concepts that fall within the umbrella term AWES, the most promising ones are crosswind AWES. Crosswind AWES can be divided into Flygen and Groundgen systems. 
Flygen systems consist of fixed-wing aircraft with wind turbines mounted aboard and connected to a ground station by a power-conducting cable. Their periodical crosswind flight (orthogonal to the wind direction) increases the apparent wind velocity, therefore increasing the withdrawn power from the wind turbines mounted onboard. 

Groundgen systems, on the other hand, convert mechanical energy to electrical energy on the ground. They comprise a tethered aircraft, usually termed a kite, which can have both rigid or flexible wings. The kite is connected by a tether that is coiled around a winch drum, coupled to a generator. Groundgen systems generate electricity in a production cycle with two phases: a power-producing reel-out (or traction) phase and a power-consuming reel-in (or retraction) phase. During the reel-out phase, the kite is controlled to fly in a fast crosswind motion, which increases the apparent wind velocity, the aerodynamic forces acting on the kite, and consequently the tension force on the tether. The tensioned tether is reeled out as it forces the generator to rotate and generate power. Eventually, the tether needs to be reeled back in. During this second phase, the kite is flown in a tension-reducing manner, and the generator, operating as a motor, coils back the tether. During the reel-in phase, some power is consumed; however, by a proper flight control, which increases tether tension during the first phase and reduces it in the second phase, the consumed power becomes a fraction of the produced power, resulting in a positive net power output.

Flight control of crosswind AWES is paramount for an efficient operation of these systems. Generally, flight control architectures follow a hierarchical structure, where a supervisor controller defines the modes of operation and the references for lower-level guidance controllers. There are different control strategies in the literature, where waypoint navigation, path-following and trajectory tracking strategies have been employed. The question of how to design a power-maximising flight controller can then be formulated in two parts: what is the optimal reference (being a waypoint, geometric path, or time-dependent trajectory) for power production, and what is a good way to track said reference. 

This work is devoted to the former, the problem of reference selection.
Reference selection involves determining spatial references, such as waypoints, geometric paths, or full trajectories, for the flight controller to track during each operational phase.
Optimal Control has been the preferred approach for generating reference signals and estimating the power production of AWES, due to its ability to model complex system dynamics and to explicitly handle state and input constraints. However, because of the system’s highly nonlinear behaviour, these formulations often lead to non-convex Optimal Control Problems (OCP) with potentially many local optima. The resulting solutions can therefore be highly dependent on the availability of a suitable initial guess, as well as on the accuracy of the system model, parameters, and disturbance representations. Furthermore, as the considered time horizon or the number of crosswind loops in the pumping cycle increases, the computational effort required to obtain a solution grows significantly.
For this reason, most of the existing literature considers periodic OCPs, in which the terminal system states are constrained to match their initial values \cite{Houska_2007, Ockels_2008, Horn_2013, Licitra_2019}. Several strategies have been proposed to facilitate the solution of OCPs and to reduce their computational burden:
\begin{itemize}
    \item \cite{Paiva_2018} use an adaptive mesh refinement scheme in a multiple-shooting direct collocation method, iteratively refining an initial coarse discretisation based on the error of the dual system.
    \item \cite{Harzer_2025} part from the assumption that crosswind loops are quasi-periodic throughout the traction phase of a PKG and adopt a stroboscopic averaging method in which they compute the result for a small subset of crosswind loops and reconstruct the full trajectory based on this subset.
    \item \cite{Gros_2013} defines a modified OCP with relaxed dynamic constraints. This relaxation is gradually tightened to match the original system dynamics. 
    \item \cite{Trevisi_2022_HBM} take the assumption of the periodicity of the solution further and describes the OCP in the frequency domain through a harmonic balance method. They then proceed to solve increasingly more complex OCPs, taking into account different physical phenomena while optimising for the mean thrust, mechanical and electrical power of a Flygen system.
    \item Finally, \cite{Erhard_2015_quaternion} proposes a singularity-free quaternion-based model OCP, thus avoiding trigonometric functions and boosting the speed of the numerical computations in the optimal control algorithm.
\end{itemize}

Other research aims to transcribe time-dependent trajectories obtained using optimal control into time-independent geometric paths, which can serve as references for path-following controllers. \cite{Fernandes_2021} proposes defining a time-independent path by extracting the dominant frequencies of the azimuth and elevation angle time profiles from the OCP solution using a Fourier transform. These frequencies are then used in the parametrisation of a Lissajous curve. \cite{Höhl_2023} compare trajectory-tracking and path-following MPC controllers, defining the reference geometric path for the path-following MPC by using cubic spline interpolation of a reference trajectory, where the time domain is resampled to a normalised interval, resulting in a parametric geometric path.

Besides OCP techniques, several research groups have been approaching the problem using adaptive or iterative learning techniques to optimise a reference path or waypoints. These techniques are able to be solved online and avoid the computationally expensive optimal control methods. 
In \cite{Zgraggen_2015}, the authors adaptively adjust the central point of a periodic eight-shaped path based on measurements of the average tether force over a full crosswind loop, thus avoiding the exact knowledge of wind conditions. 
In \cite{Cobb_2020}, the authors alter the width and height of a path defined as a Lemniscate of Gerono, maximising an objective function that, on the one hand, maximises power production and on the other hand penalises paths with sharp turns, so that the paths remain physically feasible.
As a continuation of \cite{Cobb_2020}, and also parametrising a path as a Lemniscate of Gerono, in \cite{Baheri_2019_bayesian}, the authors employ a Gaussian Process to predict the generated average power over a loop and adjust the discretised path accordingly. In the latter work, the authors manage to demonstrate that using Bayesian optimisation, the path parameters converge faster than when employing an iterative learning control method.

This work proposes a computationally tractable method of optimising a reference geometric curve, considering that the main objective of flight control is to guide the aircraft in a power-maximising crosswind motion. The path-planning problem that is analysed focuses on the reel-out phase of a Groundgen system, but a similar formulation can easily be followed for Flygen systems.
The proposed framework formulates the path-planning problem as a Nonlinear Program (NLP) which optimises the average power production during a crosswind loop over a Lissajous curve. The algorithm uses an analytic physical model of the system, based on the one presented in \cite{Schmehl_2013}, and includes a maximum curvature constraint to ensure a feasible reference path. By restricting the path to a Lissajous curve, we reduce the number of decision variables in the optimisation problem to the curve parameters, and by using a simplified physical model, we are able to guarantee a closed-form solution for the estimated power equation, therefore guaranteeing the computational efficiency of the proposed approach.

Section \ref{sec:analytic_model} describes the model used to describe the physics of a kite flying crosswind and the path parametrisation. Section \ref{sec:NLP} presents the considered optimisation problem, the applied method and results over different tether lengths. Section \ref{sec:pp_considerations} develops some considerations about the model and methodology employed, and Section \ref{sec:3_conclusions} draws concluding remarks.

\section{Analytic Physical Model}\label{sec:analytic_model}
Analytic models for AWES offer a computationally efficient manner of estimating the generated power and identifying the main parameters that influence the system's efficiency. They vary significantly in the physical effects and simplifications that are considered, but generally neglect inertial phenomena, thus usually overestimating the power potential of the system. Nevertheless, they constitute a simple, physically informed tool for the initial phases of system and operational design.

The most referenced and foundational work that laid the groundwork for subsequent analytic models is Loyd's \cite{Loyd_1980} crosswind kite power analysis for both Flygen and Groundgen systems (Drag and Lift Power Mode in Loyd's terminology). Loyd considered an idealised massless system, unencumbered by tether drag or weight, flying perfectly crosswind and neglecting any losses in the power conversion mechanism. Under these conditions, the theoretical maximum traction power that can be extracted by a Groundgen system is
\begin{equation}\label{eq:loyd}
    P_{Loyd}=\frac{1}{2}\rho A c_L v_w^3 \frac{4}{27} \left(\frac{c_L}{c_D}\right)^2,
\end{equation}
where $\rho$ is the air density, $A$ is the wing area of the kite, $c_L$ and $c_D$ are the kite's lift and drag coefficients, respectively, and $v_w$ is the wind speed.

Later works have expanded this model based on the equilibrium of forces at the kite by considering the effects of a non-perfect crosswind flight due to elevation and azimuth deviations (cosine losses), the kite and tether mass, added drag due to the tether, and considering the power consumption during the retraction phase, thus computing the net power production over a cycle of a PKG \cite{Argatov_2009, Schmehl_2013, Luchsinger_2013, Vandervlugt_2019, Trevisi_2020_unified}.  

More recently, alternative approaches have used induction and vortex models to estimate the power production of a kite, therefore considering the effects that a flying kite has on the air flow and guaranteeing adherence to the conservation of energy laws. 
\cite{delellis_2018} adapts the linear momentum models used to derive the Betz limit in conventional wind turbines to AWES. 
\cite{Trevisi_2023_refining} uses a vortex wake model developed in previous research by the same authors to refine a power equation, which can be used to optimise wing geometry and to study the performance of different systems. 
Besides power estimation and system design, induction models can also be used for kite farm layout optimisation considering the wake effects between units, while avoiding other computationally demanding methods, such as computational fluid dynamics \cite{Kheiri_2022, Karakouzian_2022}.

To obtain a closed-form expression for the power equation and prioritising analytic tractability, the model used for path-planning builds on the formulation proposed by \cite{Schmehl_2013}, which has also been used, with small variations, for other optimisation problems, such as in \cite{Roque_2020}, \cite{Rui_2025}, and more recently in \cite{da_costa_airborne_2026} to optimise the layout of an AWES farm, and in \cite{Joshi_2024} to estimate the power curve of a PKG. 
This model considers the following ubiquitous assumptions in low-fidelity physical analytic models of AWES:
\begin{itemize}
    \item Negligible gravity and inertial forces;
    \item Taut, massless tether;
    \item Constant horizontal wind;
    \item Constant kite lift and drag coefficients ($c_L$ and $c_D$);
    \item All forces acting on the kite's centre of mass.
\end{itemize}

Including inertial and gravitational forces would break the geometric similarity between aerodynamic forces and apparent wind velocity components (explained in Appendix \ref{sec:geometric_similarity}), which is the basis of a closed-form solution for the generated power equation \cite{Schmehl_2013}. This simplification is arguably acceptable for systems with negligible effects of mass, such as flexible-wing or small-scale fixed-wing systems, but may limit the direct applicability of this physical model to large-scale fixed-wing systems. However, a similar path-planning framework could be used with modifications as an iterative approach to include the effect of gravity, which is also proposed in \cite{Schmehl_2013}. 
Although gravitational and inertial forces are neglected in the aerodynamic formulation, the mass of the kite is still considered in a maximum curvature constraint, which imposes turning rate limits to the reference path.

The assumption of a taut, massless tether avoids the need for modelling tether sag and catenary effects. This is generally acceptable during the high-tension traction phase, but becomes less accurate during retraction or transition phases.  
Assuming a constant horizontal wind neglects the effect of gusts, turbulence, wind shear, and the vertical components of the wind velocity. The remaining assumptions regarding all forces acting on the centre of mass of the kite and constant lift and drag coefficients imply a point-mass model of the kite, where rotational moments are ignored, and the angle of attack is considered uniform both along the span of the aircraft and throughout the crosswind loop.

\subsection{Coordinate Systems and Path Parametrisation}

\begin{figure}[h]
    \centering
    \includegraphics[width=0.5\linewidth]{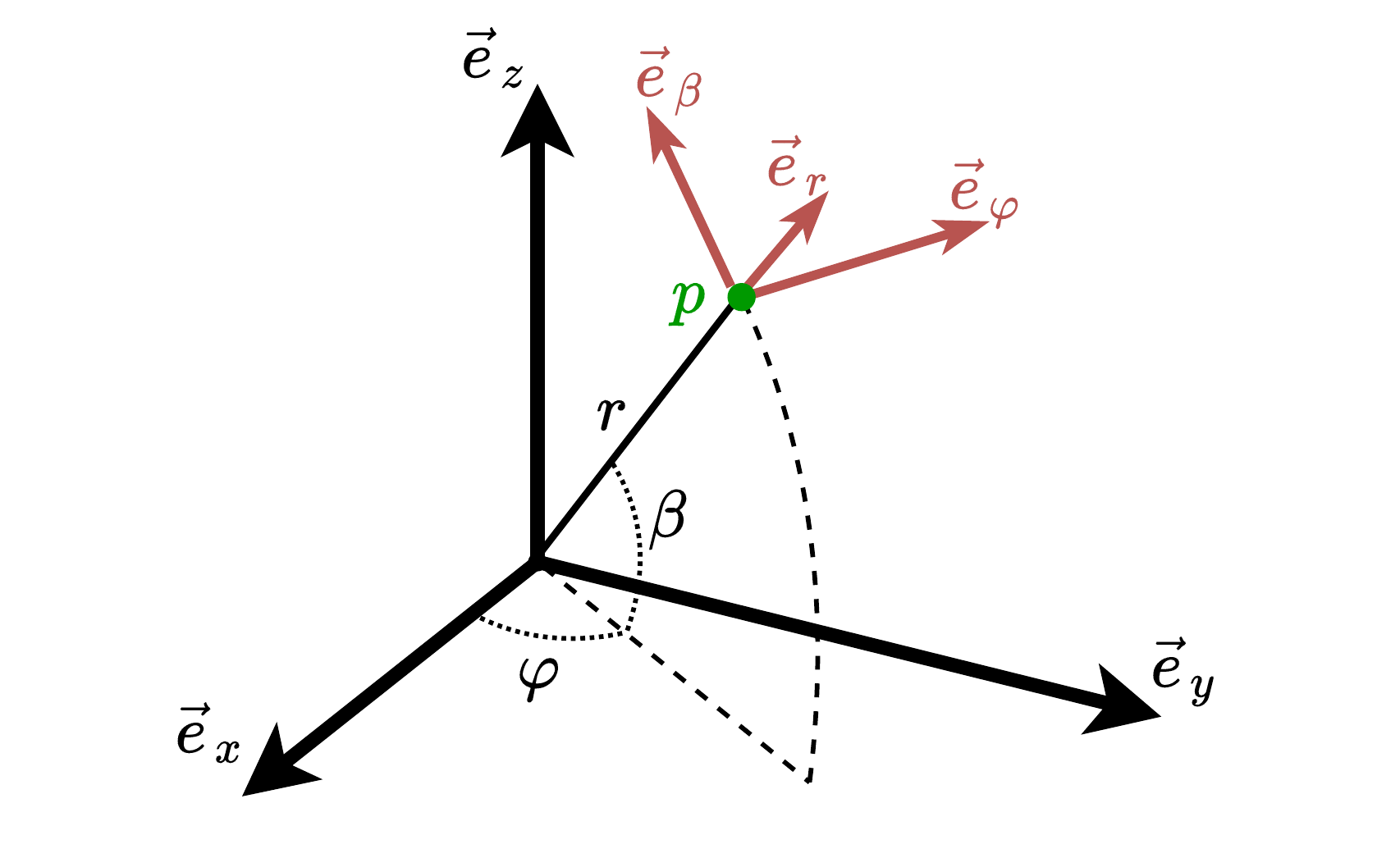}
    \caption{Global and Local coordinate systems.}
    \label{fig:global_local}
\end{figure}

The model is expressed in the following two coordinate systems:
\begin{enumerate}
    \item \textbf{Global G}: An inertial cartesian coordinate system $(x,y,z)$, with basis $(\vec{e}_x,\vec{e}_y,\vec{e}_z)$ centred at the ground station and where $z$ points vertically upwards and the $x$ axis is aligned with the predominant horizontal wind direction $\vec{v}_w=(v_w,0,0)$.
    \item \textbf{Local L}: A non-inertial spherical coordinate system $(r,\varphi,\beta)$ with basis $(\vec{e}_r,\vec{e}_\varphi,\vec{e}_\beta)$, where $r$ is the tether length (considered perfectly taut), $\varphi$ is the azimuth angle relative to the $x$ axis, and $\beta$ is the elevation angle.
\end{enumerate}
These coordinate systems are depicted in Fig. \ref{fig:global_local}. 

Considering that the kite must fly in a downwind hemisphere with positive elevation angle and the tether is perfectly taut, we can consider that the kite is free to fly over a spherical quadrant, therefore the Local coordinate system in spherical coordinates comes as the natural choice to describe the intended kite motion and the intended reference path.

\begin{figure}[h]
    \centering
    \includegraphics[width=0.7\linewidth]{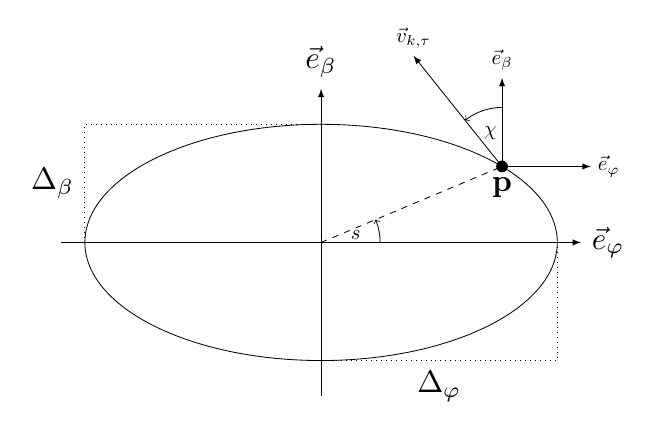}
    \caption{Path parametrisation in the $(\varphi,\beta)$ plane}
    \label{fig:path_param}
\end{figure}

This geometric path will serve as a reference for periodic crosswind loops. We have chosen a Lissajous curve parametrisation in a $(\varphi,\beta)$ plane, which allows us flexibility in terms of path design choice (the most pursued being ellipses or figures-of-eight) while reducing the amount of optimisation variables to just three parameters corresponding to the path's central elevation $\beta_0$ and its elevation and azimuthal ranges, $\Delta_\beta$ and $\Delta_\varphi$, respectively. 
The geometric path can therefore be defined as
 \begin{equation}\label{eq:path_param_chp3}
\left\{
\begin{aligned}
\beta(s) & =\beta_0 + \Delta_\beta \sin\left(\frac{n_\beta}{n_\varphi}s\right) \\
\varphi (s)  & = \Delta_\varphi \cos(s)\\
s&\in[0,2\pi]
\end{aligned}
\right.
\end{equation}
where $s$ is an angular parameter spanning between $0$ and $2\pi$, $n_\beta$ and $n_\varphi$ are the number of cycles the path describes over each axis, and its ratio $\frac{n_\beta}{n_\varphi}$ is the relative cycle rate along both coordinates. This ratio defines the number of lobes in the geometric path. A closed geometric path is only obtained for rational ratios. As an example, an ellipse would have $\frac{n_\beta}{n_\varphi}=1$ while a lying eight would have $\frac{n_\beta}{n_\varphi}=2$.
This parametrisation is shown for an elliptical path in the $(\varphi,\beta)$ plane in Fig. \ref{fig:path_param}.

\subsection{Acting Forces}

The system is considered to be dominated by the aerodynamic ($\vec{F}_{aer}$) and tether force ($\vec{F}_{tether}$). The assumption of a perfectly taut tether causes the tether force to be radial. The aerodynamic force is divided into the lift and drag forces as
\begin{equation}
    \vec{F}_{aer}=\vec{F}_{lift}+\vec{F}_{drag}
\end{equation}
where
\begin{align}
    F_{lift}&=\frac{1}{2}\rho A c_L v_a^2\\
    F_{drag}&=\frac{1}{2}\rho A c_D v_a^2
\end{align}
where $v_a$ consists of the apparent wind speed, the relative wind speed at the kite. The apparent wind velocity can be computed as
    \begin{equation}
        \vec{v}_a=\vec{v}_w-\vec{v}_k
    \end{equation}
    where $\vec{v}_w$  and $\vec{v}_k$ are the wind and the kite velocity vectors. 
    The drag force acts in the direction of $\vec{v}_a$, while the lift force acts perpendicularly to this direction.

    In \cite{Schmehl_2013}, the aerodynamic force is considered to balance the tether force perfectly and therefore
    \begin{equation}
        \vec{F}_{aer}=-\vec{F}_{tether}.
    \end{equation}
    This force balance is the basis for the derivation of the closed-form power model. 
    In this work, we divide the resulting aerodynamic force into a power-producing radial component ($\vec{F}_{aer,r}$) and a tangential steering component ($\vec{F}_{aer,\tau}$) that is responsible for the kite following the intended path. 
    While the radial component of the aerodynamic force balances the tether force and is used to derive a similar model to the one in \cite{Schmehl_2013}, a component caused by a rotation around the longitudinal axis of the kite creates a turning lift, responsible for conveying a lateral acceleration and steering the kite to follow the path.

    The radial aerodynamic force component that balances the tether force is then computed as
    \begin{equation}\label{eq:force_equilibrium}
        \vec{F}_{aer,r}=\vec{F}_{lift} \cos{(\phi)} + \vec{F}_{drag}=-\vec{F}_{tether}
    \end{equation}
    where $\phi$ is the roll angle, measuring the rotation around the kite's longitudinal axis. We must note that it is the $\vec{F}_{aer,r}$ that is fully radial, while its components, $\vec{F}_{lift}\cos{(\phi)}$ and $\vec{F}_{drag}$, may have tangential components which cancel out. 
    
\begin{figure}[h]
    \centering
    \includegraphics[width=0.6\linewidth, angle=-10]{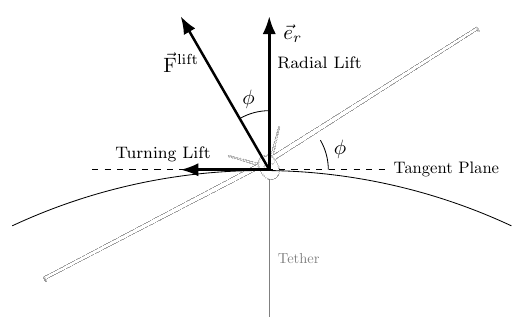}
    \vspace{-1cm}
    \caption{Turning lift on a spherical surface.}
    \label{fig:roll_angle}
\end{figure}

    The steering component of the aerodynamic force, termed the turning lift and portrayed in Fig. \ref{fig:roll_angle}, is computed as
    \begin{equation}
        \vec{F}_{aer,\tau}=\vec{F}_{lift}\sin{(\phi)}.
    \end{equation}
    The lateral acceleration that is caused by the turning lift is computed as
    \begin{equation}\label{eq:roll_acceleration}
        a_l=\frac{F_{lift}}{m}\sin{(\phi)}.
    \end{equation}

We can also compute the required lateral acceleration for the kite to follow the path as
\begin{equation}\label{eq:lateral_acceleratio}
    a_l=\frac{v_k^2}{R_0(s)}=v_k^2 \kappa(s)
\end{equation}
where $R_0$ is the local curvature radius of the path and $\kappa$ is the curvature of the path $\left(\kappa(s)=\frac{1}{R_0(s)}\right)$. The computation of $\kappa(s)$ for the Lissajous path in Eq. \eqref{eq:path_param_chp3} follows the geodesic curvature computation of Eq. \eqref{eq:curvature_geodesic} derived in the Appendix.
The roll angle that is required along the path is computed by equating Eq. \eqref{eq:roll_acceleration} and \eqref{eq:lateral_acceleratio} and solving for $\phi$, as
\begin{equation}
    \phi(s)=\sin^{-1}{\left(\frac{m v_k^2 \kappa(s)}{F_{lift}}\right)}=\sin^{-1}{\left(\frac{m v_k^2 \kappa(s)}{\frac{1}{2}\rho A c_L v_a^2}\right)}.
\end{equation}
For simplicity's sake, to guarantee a closed-form solution by avoiding a recurrence due to the computation of the apparent wind speed (which will be established below in Eq. \eqref{eq:apparent}), here we consider that $v_k \approx v_a$, an assumption also considered in, e.g., \cite{Luchsinger_2013}, to obtain
\begin{equation}\label{eq:roll_angle}
    \phi(s)=\sin^{-1}{\left(\frac{m \kappa(s)}{\frac{1}{2}\rho A c_L}\right)}.
\end{equation}

\subsection{Apparent and Kite Velocity Decomposition}
In spherical coordinates, the kite velocity can be decomposed into a radial and a tangential component, i.e., a component in the direction of the tether and a component in a plane tangent to the spherical surface.
This decomposition yields
\begin{equation}
    \vec{v}_k=\vec{v}_{k,r}+\vec{v}_{k,\tau}
\end{equation}
    where $\vec{v}_k$ is the kite velocity and $\vec{v}_{k,r}$ and $\vec{v}_{k,\tau}$ its radial and tangential components. The direction of its tangential velocity can be described by an angle $\chi$, measured relative to $\vec{e}_\beta$, which, considering the kite moves along the path, can be computed as $\chi=\tan^{-1}\left(\frac{\frac{d\varphi(s)}{ds}}{\frac{d\beta(s)}{ds}}\right)$. The Lissajous curve parametrisation permits very simple derivations of $\frac{d\beta(s)}{ds}$ and $\frac{d\varphi(s)}{ds}$, which are presented in Appendix \ref{sec:curvature}.

The above decomposition of the kite velocity in its radial and tangential components allows us to describe the apparent wind velocity as
\begin{equation}
    \vec{v}_a=\vec{v}_w-\vec{v}_{k,r}-\vec{v}_{k,\tau}
\end{equation}
which can then be described as
\begin{equation}
    \vec{v}_a=
    \begin{bmatrix}
        \cos{(\beta)}\cos{(\varphi)}\\
        \sin{(\beta)}\cos{(\varphi)}\\
        -\sin{(\varphi)}
    \end{bmatrix}v_w-
    \begin{bmatrix}
        1\\0\\0
    \end{bmatrix}v_{k,r}-
    \begin{bmatrix}
        0\\ -\cos{(\chi)}\\ -\sin{(\chi)}
    \end{bmatrix}v_{k,\tau}.
\end{equation}

By introducing the radial and tangential kite speeds as factors of the wind speed as
\begin{align}
    f&=\frac{v_{k,r}}{v_w}\\
    \lambda&=\frac{v_{k,\tau}}{v_w}
\end{align}
We can describe the apparent wind velocity as a function of the wind speed as
\begin{equation}\label{eq:va_vector_vs_vw}
    \vec{v}_a=\begin{bmatrix}
         \cos{(\beta)}\cos{(\varphi)}-f\\
        \sin{(\beta)}\cos{(\varphi)}+\lambda \cos{(\chi)}\\
        -\sin{(\varphi)}+\lambda \sin{(\chi)}
    \end{bmatrix}v_w.
\end{equation}

Similarly to the kite velocity, we can decompose the apparent wind velocity in its radial and tangential components as
\begin{equation}\label{eq:similarity}
    \vec{v}_a=\vec{v}_{a,r}+\vec{v}_{a,\tau}.
\end{equation}
Since we disregard any gravitational or inertial forces, we can establish the following relation between the radial and tangential apparent wind velocity components,
\begin{equation}\label{eq:geometric_similarity}
    \frac{v_{a,\tau}}{v_{a,r}}=\frac{F_{lift}\cos{(\phi)}}{F_{drag}}=\frac{c_L\cos{(\phi)}}{c_D}.
\end{equation}
This relation stems from the geometrical similarity of these velocity components with the $\vec{F}_{lift}\cos{(\phi)}$ and $\vec{F}_{drag}$ vectors.
This geometric similarity argument is further detailed in Appendix \ref{sec:geometric_similarity}.

From Equation \eqref{eq:va_vector_vs_vw} we find that the radial component of the apparent wind velocity is
\begin{equation}\label{eq:va_r}
    v_{a,r}=v_w(\cos{(\beta)}\cos{(\varphi)}-f)
\end{equation}
and from the geometric similarity condition, we get
\begin{align}\label{eq:apparent}
    v_a&=\sqrt{v_{a,\tau}^2+v_{a,r}^2}=v_{a,r}\sqrt{1+\left(\frac{c_L\cos{(\phi)}}{c_D}\right)^2}\\&=v_w(\cos{(\beta)}\cos{(\varphi)}-f) \sqrt{1+\left(\frac{c_L\cos{(\phi)}}{c_D}\right)^2}.
\end{align}

The tangential apparent wind velocity component follows from Eq. \eqref{eq:va_vector_vs_vw} as
\begin{equation}\label{eq:aux_eq_1}
    v_{a,\tau}=v_w \sqrt{(\sin{(\beta)}\cos{(\varphi)}+\lambda \cos{(\chi)})^2+(-\sin{(\varphi)}+\lambda \sin{(\chi)})^2}
\end{equation}
or alternatively, from Eqs. \eqref{eq:geometric_similarity} and \eqref{eq:va_r} as
\begin{equation}\label{eq:aux_eq_2}
    v_{a,\tau}=v_w(\cos{(\beta)}\cos{(\varphi)}-f)\frac{c_L\cos{(\phi)}}{c_D}.
\end{equation}
Combining these two definitions of the tangential apparent wind speed we can solve for $\lambda$ resulting in
\begin{equation}\label{eq:aux_eq_3}
    \lambda=a+\sqrt{a^2+b^2-1+\left(\frac{c_L\cos{(\phi)}}{c_D}\right)^2(b-f)^2}
\end{equation}
where
\begin{align}
    a&=-\sin{(\beta)}\cos{(\varphi)}\cos{(\chi)}+\sin{(\varphi)}\sin{(\chi)},\label{eq:a_aux}\\
    b&=\cos{(\beta)}\cos{(\varphi)},\label{eq:b_aux}
\end{align}
as is derived in Appendix \ref{sec:tangent_speed_ratio}.
This equation, without considering the effect of the roll angle ($\phi$), has been used to derive an upper bound for physically possible azimuth and elevation angles in \cite{Schmehl_2013}, since the tangential kite velocity factor $\lambda$ must not be negative.

At the point of maximum tension ($\beta=\varphi=\phi=0$), $a=0$ and $b=1$, yielding a simplified expression of $\lambda=\frac{c_L}{c_D}(b-f)$ which can be used to estimate the maximum kite speed.

\subsection{Traction Power}
Due to the assumed force equilibrium stated in Eq. \eqref{eq:force_equilibrium}, the resulting tether force is
\begin{equation}\label{eq:tetherforce}
    F_{tether}=F_{aer,r}=\frac{1}{2}\rho A  v_a^2\sqrt{(c_L\cos{(\phi)})^2+c_D^2}.
\end{equation}
Defining a resultant aerodynamic coefficient $c_R=\sqrt{(c_L\cos{(\phi)})^2+c_D^2}$ and substituting the apparent wind speed from Eq. \eqref{eq:apparent} into Eq. \eqref{eq:tetherforce} we get
\begin{equation}\label{eq:tether_force}
    F_{tether}=\frac{1}{2}\rho A c_R \left[1+\left(\frac{c_L\cos{(\phi)}}{c_D}\right)^2\right](\cos{(\beta)}\cos{(\varphi)}-f)^2v_w^2.
\end{equation}

The generated traction power corresponds to the product of the tether force and the reel-out speed, which corresponds to the radial kite speed under a perfectly taut tether consideration, yielding
\begin{equation}\label{eq:power_vs_f}
    P=F_{tether} v_{k,r}=\frac{1}{2}\rho A c_R \left[1+\left(\frac{c_L\cos{(\phi)}}{c_D}\right)^2\right]f(\cos{(\beta)}\cos{(\varphi)}-f)^2v_w^3.
\end{equation}
By differentiating this equation with respect to $f$ and equating it to zero, we can find an optimal reeling factor $f_{opt}$, as demonstrated in Appendix \ref{sec:reelout_speed_opti}, which yields 
\begin{equation}\label{eq:f_opt}
    f_{opt}=\frac{1}{3}\cos{(\beta)}\cos{(\varphi)}.
\end{equation}

For the optimal reeling factor, the maximum generated instantaneous power is
\begin{equation}\label{eq:optimal_instantaneous_power}
    P_{opt}=\frac{1}{2}\rho A c_R \left[1+\left(\frac{c_L\cos{(\phi)}}{c_D}\right)^2\right]\left(\frac{4}{27}\cos^3{(\beta)}\cos^3{(\varphi)}\right)v_w^3.
\end{equation}
Which, approximating $c_R$ by the lift coefficient $c_L$ and for a perfectly crosswind kite ($\varphi=\beta=0$), corresponds to Loyd's maximum theoretical power defined in Eq. \eqref{eq:loyd}.
We can see as well that the misalignment caused by a kite flying at an elevation and not centred, causes a penalty proportional to the cube of the cosine of the elevation and azimuth angles, termed the cosine losses.

Following the path parametrisation in Eq. \eqref{eq:path_param_chp3}, we can define the optimal power as a function of the path parameter $s$ and compute the average traction power over a closed geometric path by computing the average power along the path as
\begin{equation}
    P_{avg}=\frac{1}{2\pi}\int_0^{2\pi}P_{opt}(s)\,ds.
\end{equation}
This function will serve as our objective function in a path-planning optimisation problem.

\section{Optimisation Problem}\label{sec:NLP}

The path-planning optimisation problem is defined as the problem of selecting the optimal Lissajous curve parameters $\{\beta_0,\Delta_\beta,\Delta_\varphi\}$ that maximise the average generated power $P_{avg}$ along the closed path, subject to curvature and elevation constraints.

This can be formulated as
\begin{align}\label{eq:nlp_power}
\min_{\{\beta_0,\Delta_\beta,\Delta_\varphi\}} \quad & P_{avg} = \frac{1}{2\pi}\int_0^{2\pi}P_{opt}(s)\,ds\\
\text{s.t.} \quad & |\kappa(s)| \le \kappa_{max} \label{eq:maximum curvature}\\
&  \beta_0+\Delta_\beta \leq \beta_{max}\label{eq:maximum_beta}\\
& \beta_0-\Delta_\beta\geq \beta_{min} \label{eq:minimum_beta}\\
& \beta_0\in [\beta_{0,min}, \beta_{0,max}]\\&\Delta_\beta \in [\Delta_{\beta,min},\Delta_{\beta,max}]\\&\Delta_\varphi \in [\Delta_{\varphi,min},\Delta_{\varphi,max}]
\end{align}

The constraint in \eqref{eq:maximum curvature} limits the curvature the path can have. The value of $\kappa_{max}$ can be computed by establishing a maximum roll angle $\phi_{max}$ and solving for the $\kappa_{max}$ using Eq. \eqref{eq:roll_angle}, yielding
\begin{equation}\label{eq:max_curvature}
     \kappa_{max}=\frac{\rho A c_L }{2 m}\sin{(\phi_{max})}.
\end{equation}

\begin{figure}[h]
    \centering
    \includegraphics[width=0.5\linewidth]{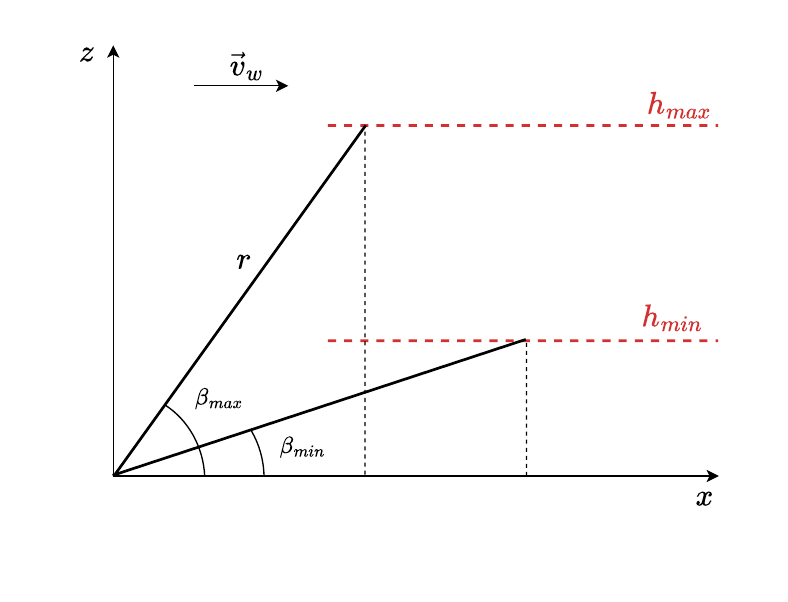}
    \caption{Maximum and minimum elevation angle constraints depiction.}
    \label{fig:height_constraint}
\end{figure}

Furthermore, we add minimum and maximum height restrictions, Eqs. \eqref{eq:maximum_beta} and \eqref{eq:minimum_beta}, to avoid possible collisions with low altitude structures or with any other aviation obstacles.  These restrictions are defined as a maximum and minimum altitude ($h_{max}$ and $h_{min}$) that are then translated to constraints on the elevation angle as portrayed in Fig. \ref{fig:height_constraint}.

The resulting constraints are
\begin{align}
    \beta_0+\Delta_\beta&\leq \beta_{max}=\sin^{-1}{\left(\frac{r}{h_{max}}\right)}\\
    \beta_0-\Delta_\beta&\geq \beta_{min}=\sin^{-1}{\left(\frac{r}{h_{min}}\right)}.
\end{align}

Finally, we add lower and upper bounds for the decision variables. These bounds serve to ensure the optimal solution is searched within the quadrant sphere corresponding to $\beta \in [0, \frac{\pi}{2}]$ and $\varphi \in [-\frac{\pi}{2}, \frac{\pi}{2}]$ and to guarantee that the found solution complies with the curvature constraint. Because, during the numerical optimisation, the curvature constraint is enforced only in a discrete sampled set of values, the optimiser can drive the decision variables $\Delta_\beta$ or $\Delta_\varphi$ to zero while satisfying the curvature constraint in the sampled set, while the continuous constraint would be violated. A small but positive minimum value $\Delta_{\beta,min}$ and $\Delta_{\varphi,min}$ guarantees a strictly positive value for these decision variables.

\subsection{Optimisation Results}
For the definition of the reference path throughout the traction phase of the PKG, we need to solve an NLP for the problem formulated in Eq. \eqref{eq:nlp_power} for several tether lengths. 
Therefore, we have set and defined an NLP for tether length increments of $5~\mathrm{m}$, starting with a $100~\mathrm{m}$ long tether until a $200~\mathrm{m}$ long tether for a small-scale kite based on the UPWIND project prototype \cite{noauthor_upwind_nodate}.

To enhance the speed of these solutions, the initial guess for each NLP was set to be the solution of the previous problem. The initial guess for the first NLP was checked a priori to guarantee its feasibility. The NLP was solved using a Sequential Quadratic Programming (SQP) algorithm using MATLAB \cite{MATLAB}.
The Lissajous curve parameter solutions $\{\beta_0, \Delta_\beta, \Delta_\varphi\}$ for each tether length were then interpolated using a cubic polynomial interpolation to establish a continuous reference path for varying tether lengths.

The maximum roll angle (used to define the curvature upper bound $\kappa_{max}$), the lower and upper bounds for the decision variables, the maximum and minimum altitude values, and other parameters are presented in Table \ref{tab:nlp_path_parameters}. 

\begin{table}[h]
\centering
\caption{Path-planning problem parameters.}
\label{tab:nlp_path_parameters}
\begin{tabular}{ll}
\hline
\textbf{Parameter} & \textbf{Value} \\
\hline
Air Density&$\rho=1.225~\mathrm{kg/m^3}$\\
Wind Speed&$v_w=10~\mathrm{m/s}$\\
Kite Mass&$m=1~\mathrm{kg}$\\
Kite Area&$A=0.28~\mathrm{m^2}$\\
Lift Coefficient&$c_L=1.2$\\
Drag Coefficient&$c_D=0.12$\\
Maximum Roll Angle & $\phi_{max} = 30^\circ$ \\
Minimum Altitude & $h_{min}=30~\mathrm{m}$\\
Maximum Altitude&$h_{max}=150~\mathrm{m}$\\

\hline
\end{tabular}
\end{table}

\begin{figure}[htbp]
    \centering
    
    \begin{subfigure}[b]{0.6\textwidth}
        \centering
        \includegraphics[width=\textwidth]{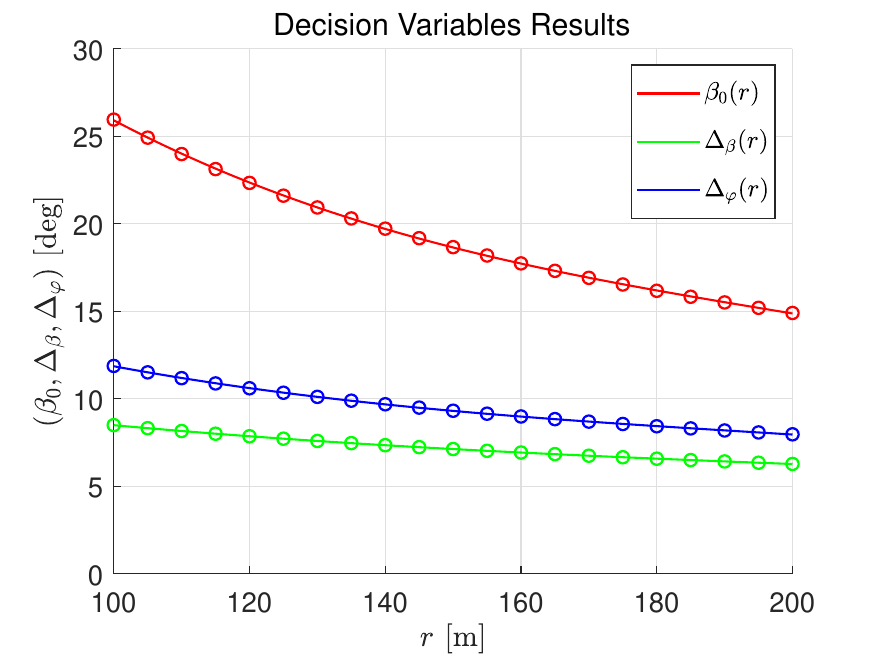}
        \caption{Decision variables results for elliptical paths.}
        \label{fig:sub1}
    \end{subfigure}
    \hfill
    \begin{subfigure}[b]{0.6\textwidth}
        \centering
        \includegraphics[width=\textwidth]{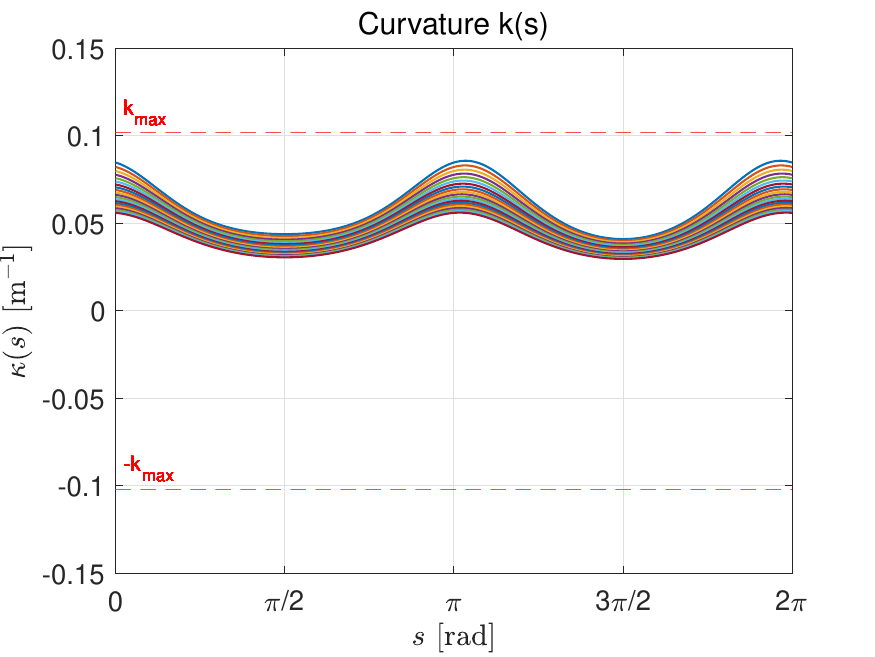}
        \caption{Curvature constraint verification for elliptical paths.}
        \label{fig:sub2}
    \end{subfigure}
    \hfill
    \begin{subfigure}[b]{0.6\textwidth}
        \centering
        \includegraphics[width=\textwidth]{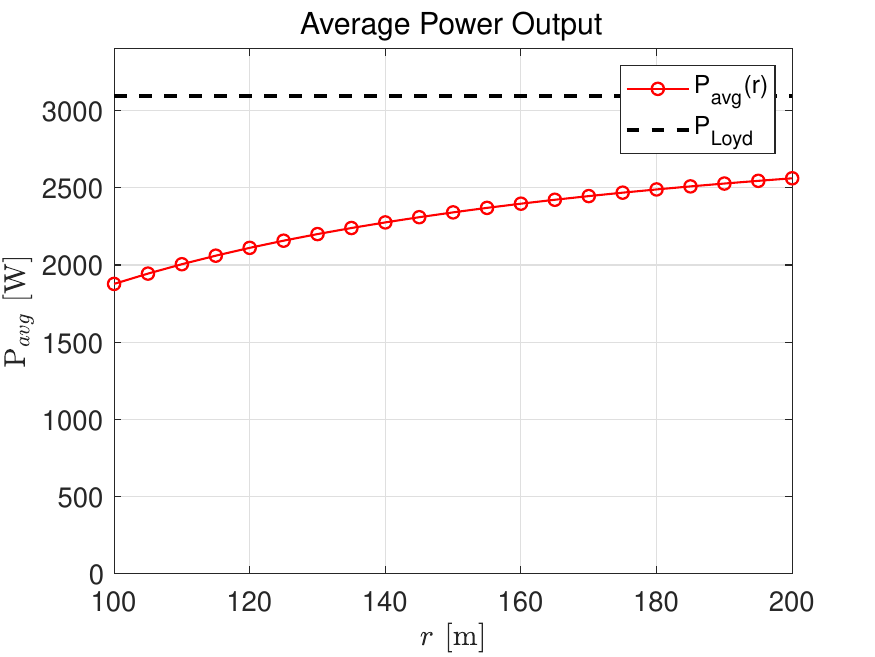}
        \caption{Average generated power with respect to tether length for elliptical paths.}
        \label{fig:sub3}
    \end{subfigure}

    \caption{Elliptical paths results.}
    \label{fig:NLP_ellipse_result}
\end{figure}

Figures \ref{fig:NLP_ellipse_result} and \ref{fig:NLP_eight_fig_result} show the resulting interpolation of $\{\beta_0, \Delta_\beta, \Delta_\varphi\}$ for tether lengths varying between $100~\mathrm{m}$ and $200~\mathrm{m}$, the curvature constraint compliance, and the estimated power compared to Loyd's theoretical limit, established in Eq. \eqref{eq:loyd}, for elliptical and figure-of-eight paths, respectively.

\begin{figure}[htbp]
    \centering
    
    \begin{subfigure}[b]{0.6\textwidth}
        \centering
        \includegraphics[width=\textwidth]{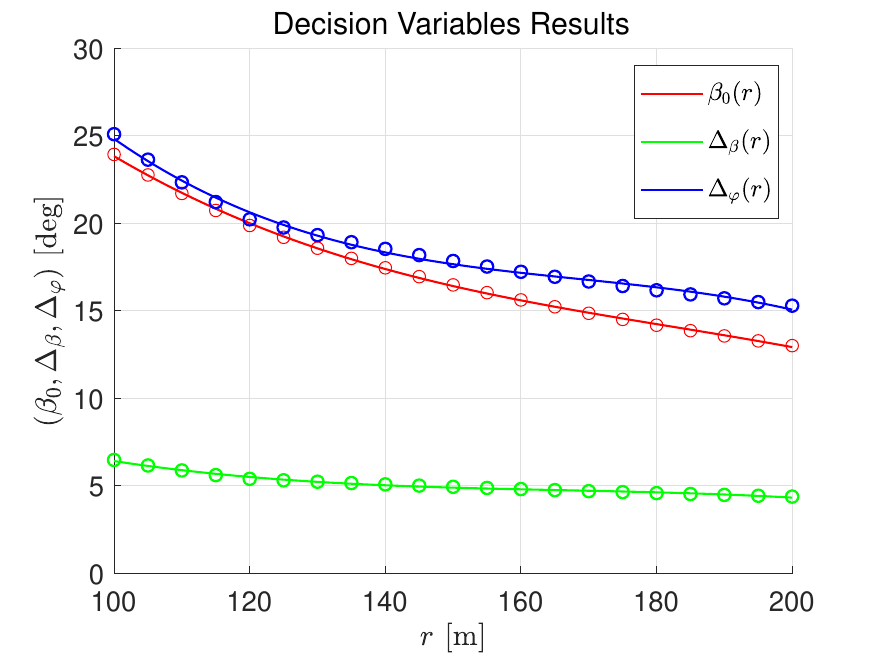}
        \caption{Decision variables results for figure-of-eight paths.}
    \end{subfigure}
    \hfill
    \begin{subfigure}[b]{0.6\textwidth}
        \centering
        \includegraphics[width=\textwidth]{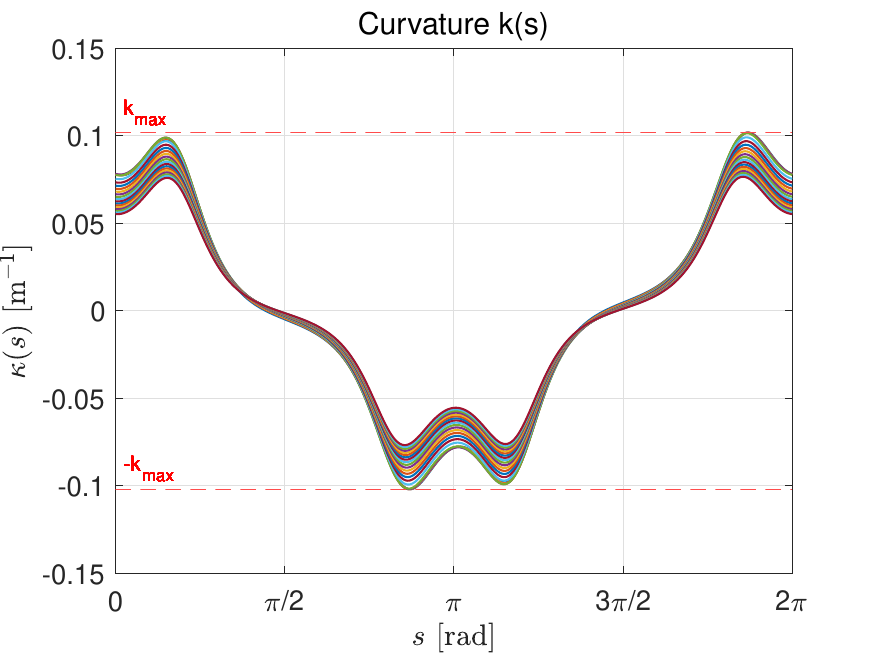}
        \caption{Curvature constraint verification for figure-of-eight paths.}
    \end{subfigure}
    \hfill
    \begin{subfigure}[b]{0.6\textwidth}
        \centering
        \includegraphics[width=\textwidth]{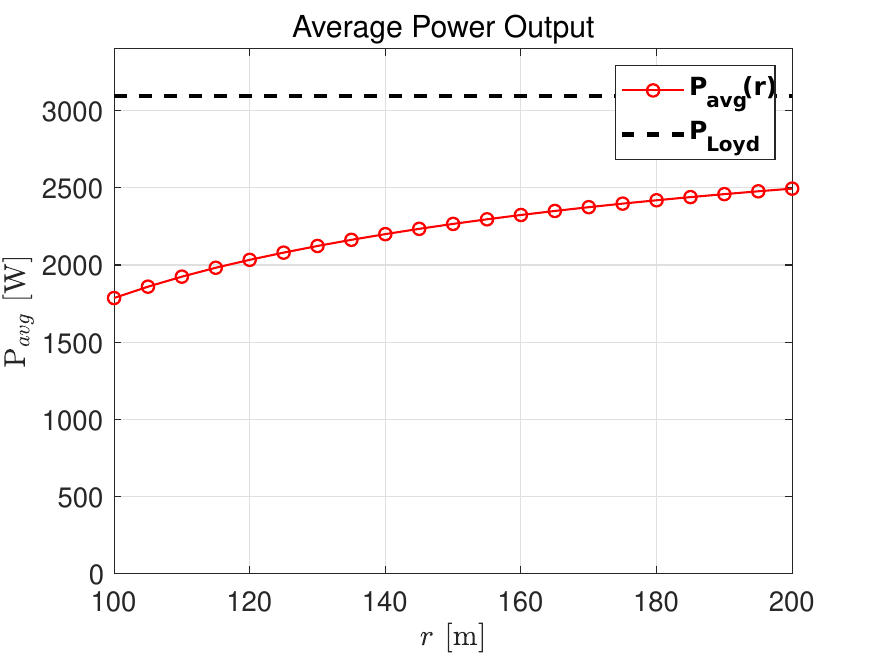}
        \caption{Average generated power with respect to tether length for figure-of-eight paths.}
    \end{subfigure}

    \caption{Figure-of-eight paths results.}
    \label{fig:NLP_eight_fig_result}
\end{figure}

\begin{figure}[htbp]
    \centering
    
    \begin{subfigure}[b]{0.7\textwidth}
        \centering
        \includegraphics[width=\textwidth]{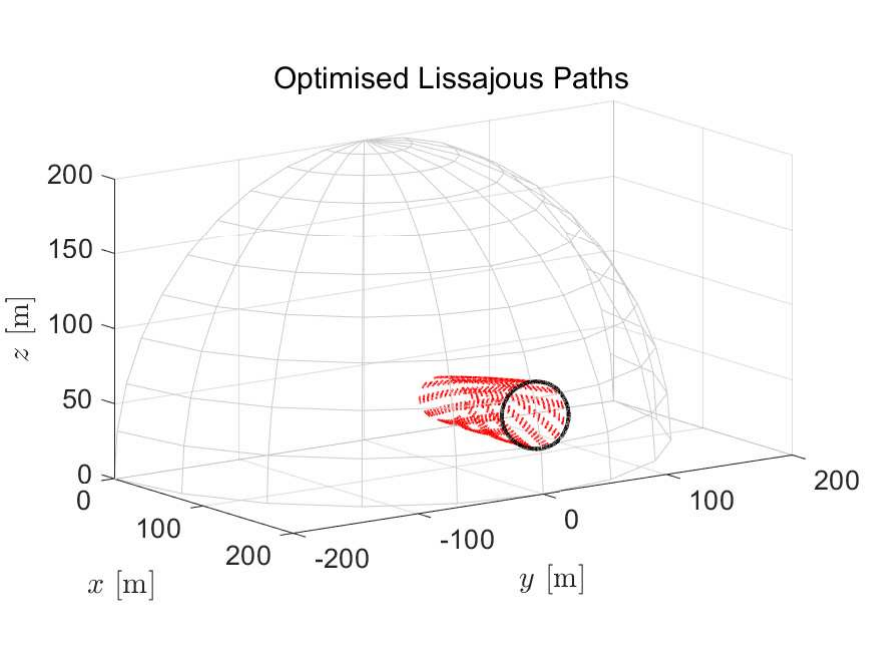}
        \caption{Elliptical optimal paths.}
    \end{subfigure}
    \hfill
    \begin{subfigure}[b]{0.7\textwidth}
        \centering
        \includegraphics[width=\textwidth]{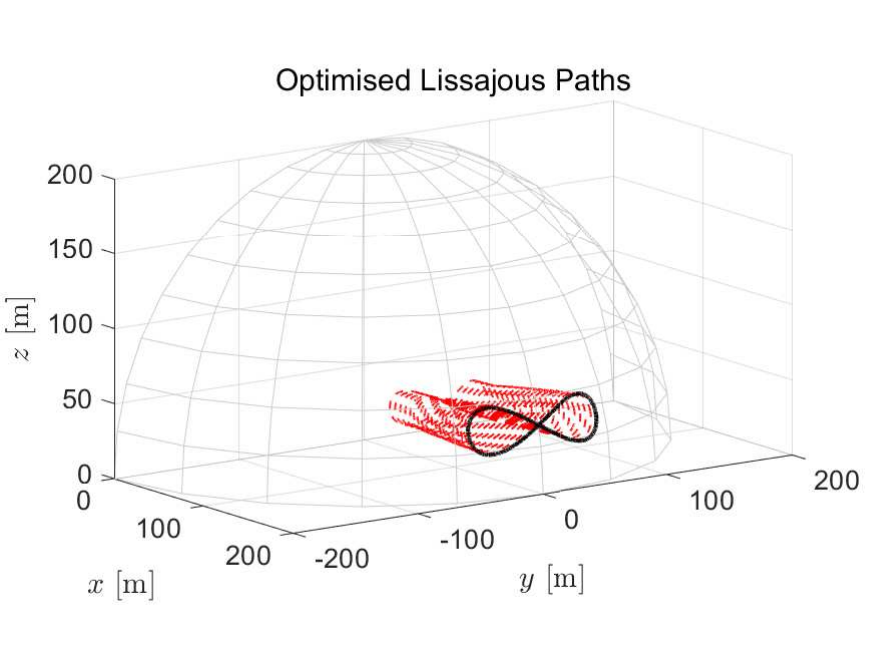}
        \caption{Figure-of-eight optimal paths.}
    \end{subfigure}

    \caption{Elliptical and figure-of-eight paths in 3D.}
    \label{fig:NLP_eight_fig_result}
\end{figure}
The results shown above reveal several consistent trends as the tether length ($r$) increases. Firstly, the elevation and azimuth angles amplitudes ($\Delta_\beta$ and $\Delta_\varphi$) decrease with longer tethers. This behaviour arises because the curvature constraint becomes less restrictive for larger values of $r$. Longer tethers reduce the geometric curvature of the path on the spherical quadrant, allowing the optimiser to choose smaller angular ranges. 

A similar trend is witnessed for the average elevation angle $\beta_0$. As the tether length increases, the average angle decreases. This is due to the minimum elevation angle constraint, which depends on a set minimum altitude. Following Eq. \eqref{eq:minimum_beta}, for larger values of $r$, the resulting $\beta_{min}$ decreases for a fixed $h_{min}$. 

These reductions of $\{\beta_0,\Delta_\beta,\Delta_\varphi\}$ for longer tethers are motivated by the considered objective function in Eq. \eqref{eq:nlp_power}, which mainly penalises deviations from a perfectly crosswind path (i.e. $(\beta,\varphi)=(0,0)$) due to the termed cosine losses. As a result, as the tether length increases, the generated power approximates Loyd's theoretical limit established in Eq. \eqref{eq:loyd}. For $r=200~\mathrm{m}$, the resulting average power reaches $83\%$ of the ideal possible generated power for elliptical paths and $81\%$ for figures-of-eight, while for $r=100~\mathrm{m}$ the generated power only reaches values of $61\%$ and $58\%$ of Loyd's limit, respectively. Further supporting this analysis, we see that the only constraint that is consistently active is the minimum elevation angle constraint. 

Furthermore, we can observe that the produced power in an elliptical path is slightly larger than the power produced for a figure-of-eight path. This is due to the larger curvature and sharper turns that are required in a figure-of-eight path that result in a larger component of the aerodynamic force responsible for steering ($\vec{F}_{aer,\tau}$) and therefore a reduction of the radial component of the aerodynamic force that contributes to power production ($\vec{F}_{aer,r}$).

The trends listed above would possibly change if other physical phenomena were to be considered. Incorporating tether drag and mass, for instance, would penalise longer tethers more strongly. These effects would add increasing aerodynamic and gravitational loads that would scale with the tether length, thus potentially affecting the monotonic increase of power production with $r$. 
On the other hand, if wind shear (the variation of wind speed with altitude) were to be added, this could result in a possibly beneficial effect on the power production with increasing tether lengths. Longer tethers allow the kite to reach higher altitudes and stronger winds. This effect would most likely benefit higher elevation angles as well, possibly counteracting the observed trend of decreasing values of $\beta_0$.
Finally, the inclusion of gravitational forces would also favour higher elevation angles so that the aerodynamic forces could better compensate for gravity, which would change the optimal trade-offs obtained in this study.

\subsection{Other Considerations}\label{sec:pp_considerations}

The curvature and roll angle computation assumes values of $v_a\approx v_k$. This assumption, although accepted in low-fidelity models, may cause possible inaccuracies. However, since $\phi$ appears in the objective function as the object of a cosine function and in the curvature computation as the object of a sine function, the possible errors caused by this simplification are further dampened. Therefore, this simplification can be deemed adequate for this problem formulation with reduced impacts on the results. 

Further constraints, such as tether force limits or rated power limits can easily be added to this framework as in \cite{Luchsinger_2013}. Using the Eqs. \eqref{eq:tether_force} and \eqref{eq:optimal_instantaneous_power} we could define such constraints as
\begin{align}
    F_{tether} &\leq F_{tether, max}\\
    P_{opt}&\leq P_{rated},
\end{align}
where $F_{tether,max}$ is the maximum traction force that the tether can sustain and $P_{rated}$ the rated power of the generator.
By incorporating these limits, it would be possible to estimate a power curve for the system and to moderate the generated power under stronger winds by altering the path. Under this setup, power moderation would take use of the cosine losses by increasing the elevation angle of the path. This possibility was not considered, since it would require the design of the overall system, namely include specificities of the tether and generator, falling outside of the scope of this paper. Furthermore, other power moderation strategies, such as the active control of the angle of attack (the angle between the wing chord and the apparent wind velocity), which influences the lift and drag coefficients ($c_L$ and $c_D$), or altering the reeling speed, diverging $f$ from its optimal value $f_{opt}$ defined in Eq. \eqref{eq:f_opt}, may be more adequate to rapidly depower the system. 

Adding a penalisation related to the tether length based on its drag or added weight could further open the possibility of defining a minimum and maximum tether length ($r_{min}$ and $r_{max}$) as new decision variables. This could structure the problem of optimising over the average power output of the full traction phase as
\begin{equation}
    P_{avg}=\frac{1}{2\pi(r_{max}-r_{min})}\int_{r_{min}}^{r_{max}}\int_0^{2\pi}P_{opt}(s)\,ds \,dr.
\end{equation}

Another variant of this problem could be to add the reeling speed ratio $f$ to the decision variables to find a constant reeling speed reference, and utilise, not $P_{opt}$, but the instantaneous power computed in Eq. \eqref{eq:power_vs_f}.

Finally, a thorough analysis regarding the adequacy of the SQP to solve such a problem is lacking. This problem is smooth, and the initial guesses are guaranteed to be feasible; however, the objective function is not fully convex, and it was not verified whether a constraint qualification is satisfied. The nonconvexity of the objective function can cause the optimiser to converge to a local optimum, and the possible linearly dependent constraints could cause issues when applying this method. These issues have not been witnessed for this choice of parameters, but further research is required to either guarantee constraint qualification of this problem setup or to establish limits to this approach. Nonetheless, MATLAB's optimiser is prepared to circumvent these issues, and these were not witnessed while obtaining the presented results.

\section{Concluding Remarks}\label{sec:3_conclusions}

This work presents a computationally efficient method to compute the optimal path for crosswind flight by parametrising the path as a Lissajous curve, thereby reducing the number of decision variables to the curve parameters.
The optimisation problem is here formulated as an NLP and solved using SQP for increasing tether lengths. The results are then interpolated to find a path reference for varying tether lengths. 

This formulation not only can be used for efficient path-planning for the traction phase of a Groundgen system but also for its power estimation. In fact, a similar formulation was used in an AWE farm layout optimisation problem in \cite{Roque_2020}, \cite{Rui_2025} and \cite{da_costa_airborne_2026}. Furthermore, minor alterations can be made to adapt this formulation to Flygen systems.

The physical model used in this formulation, although allowing a closed-form solution of the power equation due to its simplicity, fails to include forces such as gravity and inertia. This may affect its applicability to large-scale fixed-wing systems. This points to the necessity of performing validation studies, comparing the results with OCP solutions, in order to analyse for which system sizes this framework is suitable. Future work may also be focused on the inclusion of these forces while retaining computational tractability. Furthermore, a sensitivity analysis will be developed to study how varying parameters and constraints can affect the shape of the obtained paths and power production. Finally, besides interpolating the path parameters over varying tether lengths, doing so for different wind speeds will allow us to achieve references that vary smoothly for different wind conditions.



\appendix
\section{Appendix - Auxiliary Equations}
This appendix serves as an aid to better comprehend the implicit derivations used in the development of the analytic physical model presented above and used for path planning. The following sections present the derivations that explain the steps to achieve the equations \eqref{eq:aux_eq_3} and \eqref{eq:f_opt}, and the geometric similarity in \eqref{eq:geometric_similarity} that permits a closed-form solution under a negligible gravity and inertia assumption. Furthermore, an explanation of how the geodesic curvature is computed is provided in Section \ref{sec:curvature}.

\subsection{Apparent Wind and Aerodynamic Forces Geometric Similarity}\label{sec:geometric_similarity}

In this model, we decompose the apparent wind velocity $\vec{v}_a$ in two orthogonal components, the radial ($\vec{v}_{a,r}$) and tangential ($\vec{v}_{a,\tau}$) components. The total aerodynamic force $\vec{F}_{aer}$ can also be decomposed into two orthogonal vectors corresponding to the force of Lift ($\vec{F}_{lift}$) and Drag ($\vec{F}_{drag}$).

By definition, the Drag force is collinear with the apparent wind velocity ($\vec{F}_{drag} ||\vec{v}_a$), and since we are neglecting any other forces besides the tether and aerodynamic forces, the force equilibrium stated in Equation \eqref{eq:force_equilibrium} guarantees that the total aerodynamic force points in the radial direction and it is collinear with the radial component of the apparent wind velocity ($\vec{F}_{aer}||\vec{v}_{a,r}$).

\begin{figure}[h]
    \centering
    \includegraphics[width=0.5\linewidth]{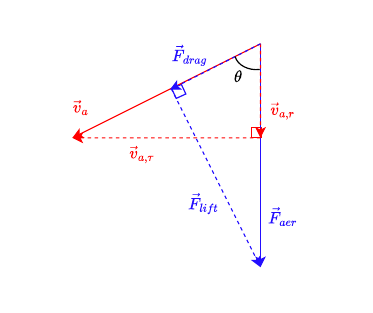}
    \caption{Geometric similarity of the aerodynamic forces and the apparent velocity components.}
    \label{fig:triangles_decomposition}
\end{figure}

Figure \ref{fig:triangles_decomposition} depicts the two right triangles formed by the apparent wind velocity and aerodynamic force decomposition rotated to lie on the same plane.

\begin{figure}[h]
    \centering
    \includegraphics[width=0.5\linewidth]{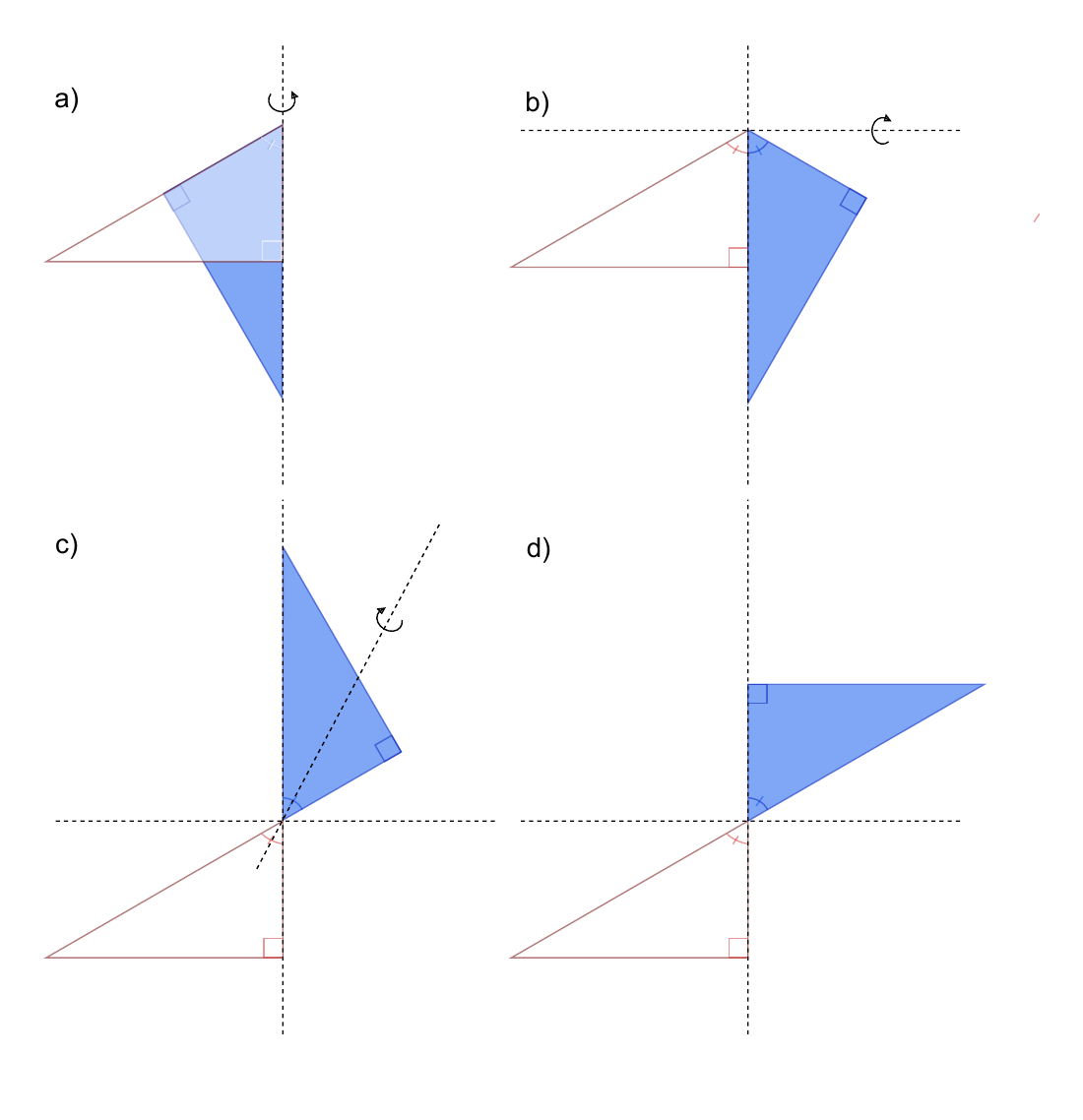}
    \caption{Geometric similarity verification by triangle rotation.}
    \label{fig:triangle_rotation}
\end{figure}

These two triangles share a common corner with an angle $\theta$ that can be computed as
\begin{equation}
    \theta=\tan^{-1}{\left(\frac{v_{a,\tau}}{v_{a,r}}\right)}=\tan^{-1}{\left(\frac{F_{lift}\cos{(\phi)}}{F_{drag}}\right)}
\end{equation}
from which we yield the resulting ratio similarity stated in Equation \eqref{eq:geometric_similarity}
\begin{equation*}
    \frac{v_{a,\tau}}{v_{a,r}}=\frac{F_{lift}\cos{(\phi)}}{F_{drag}}.
\end{equation*}
This geometric similarity can more easily be verified by fixing one of the triangles and rotating the other around some well-chosen axes, as shown in Fig. \ref{fig:triangle_rotation}.

\subsection{Tangential Kite Speed Ratio}\label{sec:tangent_speed_ratio}

The tangential kite speed ratio $\lambda$, defined in Eq. \eqref{eq:aux_eq_3}, can be computed by parting from Eqs. \eqref{eq:aux_eq_1} and \eqref{eq:aux_eq_2}. Here repeated for readability,
\begin{align*}
    v_{a,\tau}&=v_w \sqrt{(\sin{(\beta)}\cos{(\varphi)}+\lambda \cos{(\chi)})^2+(-\sin{(\varphi)}+\lambda \sin{(\chi)})^2}\\
    v_{a,\tau}&=v_w(\cos{(\beta)}\cos{(\varphi)}-f)\frac{c_L\cos{(\phi)}}{c_D}.
\end{align*}

Let us begin by defining the auxiliary variables, shown in Eq. \eqref{eq:a_aux} and \eqref{eq:b_aux},
\begin{align*}
    a&=-\sin{(\beta)}\cos{(\varphi)}\cos{(\chi)}+\sin{(\varphi)}\sin{(\chi)},\\
    b&=\cos{(\beta)}\cos{(\varphi)}.
\end{align*}

Equating \eqref{eq:aux_eq_1} and \eqref{eq:aux_eq_2} and squaring both sides of the equation, we get

\begin{equation}
    (b-f)^2 \left(\frac{c_L \cos{(\phi)}}{c_D}\right)^2 =(\sin{(\beta)}\cos{(\varphi)}+\lambda \cos{(\chi)})^2+(-\sin{(\varphi)}+\lambda \sin{(\chi)})^2,
\end{equation}
which can then be simplified as
\begin{equation}
    \lambda^2-2\lambda a+\sin^2{(\beta)}\cos^2{(\varphi)}+\sin^2{(\varphi)}-(b-f)^2 \left(\frac{c_L \cos{(\phi)}}{c_D}\right)^2=0.
\end{equation}

Since the tangential speed ratio cannot be negative, this quadratic equation yields the following solution
\begin{equation}
    \lambda=a+\sqrt{a^2+(b-f)^2 \left(\frac{c_L \cos{(\phi)}}{c_D}\right)^2-\sin^2{(\beta)}\cos^2{(\varphi)}-\sin^2{(\varphi)}}.
\end{equation}

The terms $-\sin^2{(\beta)}\cos^2{(\varphi)}-\sin^2{(\varphi)}$ can be simplified in the following manner
\begin{align}
    -\sin^2{(\beta)}\cos^2{(\varphi)}-\sin^2{(\varphi)}&=-\sin^2{(\beta)}\cos^2{(\varphi)}-\sin^2{(\varphi)}+\cos^2{(\beta)}\cos^2{(\varphi)}-\cos^2{(\beta)}\cos^2{(\varphi)}\\
    &=-\cos^2{(\varphi)}-\sin^2{(\varphi)}+\cos^2{(\beta)}\cos^2{(\varphi)}\\
    &=-1+\cos^2{(\beta)}\cos^2{(\varphi)}\\
    &=b^2-1,
\end{align}
finally yielding Eq. \eqref{eq:aux_eq_3},
\begin{equation*}
    \lambda=a+\sqrt{a^2+b^2-1+(b-f)^2 \left(\frac{c_L \cos{(\phi)}}{c_D}\right)^2}.
\end{equation*}

\subsection{Optimal Reel-out Speed Ratio}\label{sec:reelout_speed_opti}

The optimal reel-out speed ratio, or optimal radial kite speed ratio, defined in Eq. \eqref{eq:f_opt}, can be found by differentiating the traction power equation \eqref{eq:power_vs_f} with respect to the speed ratio $f$, equating it to zero and solving for $f$, thereby finding the critical points.

The derivative of \eqref{eq:power_vs_f} with respect to $f$ is
\begin{equation}
    \frac{dP}{df}=\frac{1}{2}\rho A c_R\left[1+\left(\frac{c_L\cos{(\phi)}}{c_D}\right)^2\right]\left[3f^2-4f\cos{(\beta)}\cos{(\varphi)}+\cos^2{(\beta)}\cos^2{(\varphi)}\right]v_w^3.
\end{equation}

By equating it to zero, one can further simplify it as
\begin{align}
    3f^2-f\cos{(\beta)}\cos{(\varphi)}+\cos^2{(\beta)}\cos^2{(\varphi)}=0 \Leftrightarrow\\
    f=\frac{2}{3}\cos{(\beta)}\cos{(\varphi)} \pm \frac{\sqrt{4\cos^2{(\beta)}\cos^2{(\varphi)}}}{6}\Leftrightarrow\\
    f=\frac{2}{3}\cos{(\beta)}\cos{(\varphi)} \pm \frac{2\cos{(\beta)}\cos{(\varphi)}}{6}\Leftrightarrow\\
    f=\frac{1}{3}\cos{(\beta)}\cos{(\varphi)}~\lor~f=\cos{(\beta)}\cos{(\varphi)}.
\end{align}

By replacing both possible solutions in \eqref{eq:power_vs_f}, since the second solution yields zero traction power, we conclude the optimal reeling speed factor is 
\begin{equation*}
    f_{opt}=\frac{1}{3}\cos{(\beta)}\cos{(\varphi)}
\end{equation*}
and the optimal instantaneous traction power, presented in Eq. \eqref{eq:optimal_instantaneous_power},
\begin{equation*}
    P_{opt}=\frac{1}{2}\rho A c_R \left[1+\left(\frac{c_L\cos{(\phi)}}{c_D}\right)^2\right]\left(\frac{4}{27}\cos^3{(\beta)}\cos^3{(\varphi)}\right)v_w^3.
\end{equation*}

\subsection{Geodesic Curvature of a Path on a Sphere}\label{sec:curvature}
Considering a path as a parametric curve in a three-dimensional space as a smooth mapping $\mathbf{p}_{ref}: [a,b]\subset\mathbb{R}\rightarrow\mathbb{R}^3, \mathbf{p}_{ref}(s)=(x(s),y(s),z(s))$, where $x,y,z:[a,b]\rightarrow\mathbb{R}$ are at least $C^2$ functions and $s$ is a general parameter along the curve, we can compute the curvature of the path as
\begin{equation}\label{eq:curvature_3D}
    \kappa(s)=\frac{\Vert \mathbf{p}'_{ref}(s) \times \mathbf{p}''_{ref}(s) \Vert}{\Vert \mathbf{p}'_{ref}(s) \Vert^3},
\end{equation}
where $\mathbf{p}'_{ref}(s)$ and $\mathbf{p}''_{ref}(s)$ are, respectively,
\begin{align}
    \mathbf{p}'_{ref}(s)&=
    \begin{bmatrix}
        x'(s)\\y'(s)\\z'(s)
    \end{bmatrix}\\
    \mathbf{p}''_{ref}(s)&=
    \begin{bmatrix}
        x''(s)\\y''(s)\\z''(s)
    \end{bmatrix}
\end{align}
and 
\begin{align}
    x'(s) &= \frac{\mathrm{d}x}{\mathrm{d}s}, \quad y'(s) = \frac{\mathrm{d}y}{\mathrm{d}s}, \quad z'(s) = \frac{\mathrm{d}z}{\mathrm{d}s},\\
x''(s) &= \frac{\mathrm{d}^2x}{\mathrm{d}s^2}, \quad y''(s) = \frac{\mathrm{d}^2y}{\mathrm{d}s^2}, \quad z''(s) = \frac{\mathrm{d}^2z}{\mathrm{d}s^2}.
\end{align}

The geodesic curvature $\kappa_{geo}$ of a path on a surface within this three-dimensional space corresponds to the component of the total curvature of the path, computed in Eq. \eqref{eq:curvature_3D}, in the tangent plane of the considered surface. In a spherical surface with radius $r$, the normal curvature (curvature of a curve in the normal direction of the surface) is equal to $\kappa_n=\frac{1}{r}$. Therefore, considering the total curvature tangential and normal decomposition, the geodesic curvature can be computed as
\begin{align}
   \kappa(s)^2 &= \kappa_{geo}(s)^2+\kappa_n^2 \Leftrightarrow\\
   \kappa_{geo}(s)&=\sqrt{\kappa(s)^2-\frac{1}{r^2}}.\label{eq:curvature_geodesic}
\end{align}

Considering the path parametrisation as a Lissajous curve in a $(\varphi,\beta)$ plane, as described in Eq. \eqref{eq:path_param_chp3}, we can translate this path to a three-dimensional path as
\begin{equation}\label{eq:path_global_3D}
    \mathbf{p}_{ref}(s)=r\begin{bmatrix}
        \cos{(\beta(s))}\cos{(\varphi(s)}\\
        \cos{(\beta(s))}\sin{(\varphi(s))}\\
        \sin{(\beta(s))}
    \end{bmatrix}.
\end{equation}

We can derive the first and second order derivatives of $\varphi(s)$, $\beta(s)$, which can be expressed as
\begin{equation}
\left\{
\begin{aligned}
\varphi'(s) & =  -\Delta_\varphi \sin{(s)}\\
\beta'(s)  & = \Delta_\beta \frac{n_\beta}{n_\varphi} \cos{\left(\frac{n_\beta}{n_\varphi}s\right)},
\end{aligned}
\right.
\end{equation}
and
\begin{equation}
\left\{
\begin{aligned}
\varphi''(s) & =  -\Delta_\varphi \cos{(s)}\\
\beta''(s)  & = -\Delta_\beta \left(\frac{n_\beta}{n_\varphi}\right)^2 \sin{\left(\frac{n_\beta}{n_\varphi} s\right)}.
\end{aligned}
\right.
\end{equation}
Considering these $\varphi$ and $\beta$ derivatives and the path definition as in Eq. \eqref{eq:path_global_3D}, the first and second order derivatives of $\mathbf{p}_{ref}$ yield
\begin{equation}
    \mathbf{p}_{ref}'(s)=r
    \begin{bmatrix}
        -\sin{(\beta)}\cos{(\varphi)}\beta'-\cos{(\beta)}\sin{(\varphi)}\varphi'\\
        -\sin{(\beta)}\sin{(\varphi)}\beta'+\cos{(\beta)}\cos{(\varphi)}\varphi'\\
        \cos{(\beta)}\beta'
    \end{bmatrix}
\end{equation}
and
\begin{equation}
    \mathbf{p}_{ref}''(s)=r
    \begin{bmatrix}
        -\cos{(\beta)}\cos{(\varphi)}[(\beta')^2+(\varphi')^2]+2\sin{(\beta)}\sin{(\varphi)}\beta'\varphi'-\sin{(\beta)}\cos{(\varphi)}\beta''-\cos{(\beta)}\sin{(\varphi)}\varphi''\\
        -\cos{(\beta)}\sin{(\varphi)}[(\beta')^2+(\varphi')^2]-2\sin{(\beta)}\cos{(\varphi)}\beta'\varphi'-\sin{(\beta)}\sin{(\varphi)}\beta''+\cos{(\beta)}\cos{(\varphi)}\varphi''\\
        -\sin{(\beta)}(\beta')^2+\cos{(\beta)}\beta''
    \end{bmatrix}
\end{equation}
where, for clarity, the explicit dependence on the path parameter $s$ is omitted.
From these equations, we can compute the total and geodesic curvatures of the path using Eqs. \eqref{eq:curvature_3D} and \eqref{eq:curvature_geodesic} above.

\bibliography{myrefs}

\end{document}